\documentclass[12pt]{amsart}

\usepackage[new]{old-arrows}
\usepackage{longtable}
\usepackage{booktabs}
\usepackage{array}
\usepackage{multirow}
\usepackage{supertabular}
\usepackage{amssymb}
\usepackage{amsmath}
\usepackage{mathrsfs}
\usepackage{titletoc}
\usepackage{graphicx}
\usepackage{latexsym}
\usepackage[alphabetic]{amsrefs}

\usepackage{tikz}
\usetikzlibrary{matrix}
\usepackage{tikz,tikz-cd}

\usepackage{amsmath}
\usepackage{hyperref}
\usepackage{amssymb}
\usepackage{latexsym}
\usepackage{amscd}
\usepackage{graphicx} 
\usepackage{amsthm}
\usepackage{mathrsfs}
\usepackage{xypic}
\usepackage{bm}
\usepackage{color}
\usepackage{ulem}

\newdimen\AAdi%
\newbox\AAbo%
\def\AAk#1#2{\s_etbox\AAbo=\hbox{#2}\AAdi=\wd\AAbo\kern#1\AAdi{}}%
\def\AAr#1#2#3{\s_etbox\AAbo=\hbox{#2}\AAdi=\ht\AAbo\raise#1\AAdi\hbox{#3}}%
\font\tenmsb=msbm10 at 12pt \font\sevenmsb=msbm7 at 8pt
\font\fivemsb=msbm5 at 6pt
\newfam\msbfam
\textfont\msbfam=\tenmsb \scriptfont\msbfam=\sevenmsb
\scriptscriptfont\msbfam=\fivemsb

\newtheorem{thm}{Theorem}[section]

\newtheorem{rem}[thm]{Remark}

\newtheorem{defi}[thm]{Definition}

\newcommand{\Section}[2]{\setcounter{equation}{0}
\allowdisplaybreaks
\section[#1]{#2}}

\def\pf#1{\frac{\partial}{\partial #1}}
\def\pd#1#2{\frac {\partial #1}{\partial #2}}

\def\a{\alpha}

\def\p#1{\partial #1}

\def\de{\delta}

\def\la{\lambda}

\def\ra{\rightarrow}

\subjclass[2010]{58E20,~53A10,~53C42.}
\begin{document}
\pagenumbering{Roman}\setcounter{page}{1}

\pagenumbering{arabic} \setcounter{page}{1}
\title[Dirichlet problem and minimal cones]{Recent progress on the Dirichlet problem for the minimal surface system and minimal cones}
\author{Yongsheng\ Zhang}
\address{
            Tongji University \& Max Planck Institute for Mathematics at Bonn}\email
 {
yongsheng.chang@gmail.com
}

\date{}

\thanks{Sponsored in part by NSFC (Grant No. 11601071), 
and a Start-up Research Fund from Tongji University.}

\begin{abstract}

This is a very brief report on recent developments on the Dirichlet problem for the minimal surface system and minimal cones in Euclidean spaces.
We shall mainly focus on two directions:

(1) 
Further systematic developments after Lawson-Osserman's paper \cite{l-o} on the Dirichlet problem for minimal graphs of high codimensions.
Aspects including non-existence, non-uniqueness and irregularity properties of solutions have been explored from different points of view.

(2) Complexities and varieties of area-minimizing cones in high codimensions. 
We shall mention interesting history and exhibit some recent results which successfully furnished new families of minimizing cones of different types.

\end{abstract}

\maketitle



\Section{Introduction}{Introduction}\label{S1}
\subsection{Plateau problem}
The problem is to consider minimal surfaces spanning a given contour.
Roughly speaking, there are two kinds depending on desired minimality.
One is the ``minimizing" setting for finding global minimizers for area functionals under various boundary or topological constraints;
while the other is the ``minimal" setting for critical points.

Stories of the problem can trace back to J.-L. Lagrange
who, in 1768, considered graphs with minimal area over some domain $D$ of $\mathbb R^2$.
A necessary condition is the Euler-Lagrange equation, for $z=z(x,y)$,
\begin{equation}\label{2dim}
(1+z_y^2)z_{xx}-2z_xz_yz_{xy}+(1+z_x^2)z_{yy}=0.
\end{equation}
%
From then on, the theory of minimal surfaces (with vanishing mean curvature)
soon launched an adventure journey.
 Lots of great mathematicians,
including Monge, J. Meusnie,  A.-M. Legendre, S. Poisson, H. Scherk, E. Catalan, O. Bonnet, H. Schwarz, S. Lie and many others,
 entered this filed and made it flourishing for more than one century.
 
During the mathematical developments, Belgian physicist J. Plateau did a good number of intriguing experiments with soap films (not merely using wires)
and in \cite{P} gained some explanation about the phenomena of stability and instability,
i.e., whether or not small deformations of the film can decrease area.
By laws of surface tension, an observable soap film bounded by a given simple closed curve is stable minimal.
Thus Plateau provided physical solutions to the question in $\mathbb R^3$ in the minimal setting,
and the problem was named after Plateau since then.

However, it took more time for rigorous mathematical arguments.
In 1930, J. Douglas \cite{d} and T. Rad\'o \cite{r} affirmatively answered the problem in $\mathbb R^3$, 
respectively, in the minimizing setting.  
General cases were subsequently studied and a big portion 
were solved due to Federer and Fleming's celebrated compactness theorems of normal currents and integral currents \cite{FF}
in expanded territories.
\subsection{Dirichlet problem for minimal graph of condimension one}

It can be seen that Plateau problem (for minimal surfaces with given simple closed boundary curves)
 is  actually beyond the scope of Lagrange's  original question.
 
 If $D$ is a bounded domain of $\mathbb R^{n+1}$ with $C^2$ boundary and $\phi:\p D\rightarrow {\mathbb R}^{1}$,
 then the Dirichlet problem for minimal surface equation is asking for solution $f: D\rightarrow {\mathbb R}^{1}$ to satisfy
 following generalization of \eqref{2dim}
 \begin{equation}\label{DP1}
(1+|\nabla f|^2)\triangle f- \sum_{i,j=1}^{n+1}f_i f_j f_{ij}=0
 \end{equation}
 and
 $f|_{\p D}=\phi$.
Hence the Dirichlet problem can be regarded as a special kind of Plateau problem, 
which searches for graph solutions
for graph boundary data.

For $n+1=2$ and convex $D$, the Dirichlet problem is solvable for any continuous boundary data, see 
\cite{r2}. 
In general situation, 
by efforts of Jenkins-Serrin \cite{j-s} and later Bombieri-De Giorgi-Miranda \cite{b-d-m},
the Dirichlet problem turns out to be well posed (i.e., having a unique solution) for any continuous boundary function if and only if 
$\p D$ is everywhere mean convex.
Moreover, if solution exists, it must be $C^\omega$ due to de Giorgi \cite{de} (also see \cite{St} and \cite{m1});
and its graph is absolutely area-minimizing (see \cite{fe}), 
which means any competitor sharing the same boundary possesses larger volume.

Dirichlet problem for minimal graphs of high codimensions will be discussed in \S \ref{S2}.



\subsection{Bernstein problem}
In his  paper \cite{B}, Bernstein showed that
every solution to \eqref{DP1} for $n=1$ or \eqref{2dim} in the entire $\mathbb R^2$ (with no boundary requirement at infinity)
has to be affine. Fleming \cite{fle} suggested a new idea for this problem which also works for $n\geq 2$ via De Giorgi's improvement \cite{de2}.
The principle states that the existence of a non-affine solution over $\mathbb R^{n+1}$ implies the existence of a non-planar area-minimizing hypercone in $\mathbb R^n$.
Almgren \cite{A} followed this line and  gained the same conclusion for $n=2$.
In \cite{S} J. Simons greatly extended the results 
by showing no non-planar stable hypercones in $\mathbb R^{n+1}$ for $n\leq 6$.
In $\mathbb R^8$, he discovered stable minimal hypercones
         \begin{equation}
          C_{k,k}
          =
          C\left(
          S^{k}\left(\sqrt{\frac{1}{2}}\right)
         \times 
         S^{k}\left(\sqrt{\frac{1}{2}}\right)
         \right)
         \subset \mathbb R^{2(k+1)}
         \ \text{ when } k\geq 3.
         \end{equation}
Here, for a set $E$ in unit sphere, the cone over $E$ is defined to be $C(E):=\{tx:x\in E,\ t\in(0,\infty)\}$.
Then he naturally raised the question whether $C_{k,k}$ for $k\geq 3$ in $\mathbb R^{2(k+1)}$, nowadays called Simons cones, are area-minimizing.
Immediately, the celebrated article \cite{b-d-g} by Bombieri-De Giorgi-Giusti confirmed that all Simons cones are area-minimizing and 
constructed a non-planar minimal graph over $\mathbb R^8$ in $\mathbb R^9$ which has $C_{3,3}\times \mathbb R$ as its tangent cone at infinity. 
As a result, the yes-no answer to the Bernstein problem got complete: 
there exist no non-planar minimal graphs over $\mathbb R^{n+1}$ in $\mathbb R^{n+2}$ when $n\leq 6$;
but there are such creatures when $n\geq 7$.

Still, lots of interesting subtle behaviors are  mysterious to us, 
such as what types of entire minimal graphs can occur? 
Right after \cite{b-d-g}, H. B. Lawson, Jr. considered equivariant plateau problems in \cite{l} and  obtained almost all homogeneous area-minimizing hypercones
(see \cite{zha}).
P. Simoes  \cite{PS1, PS2} added that $C_{2,4}$ is also minimizing.
R. Hardt and L. Simon \cite{HS} discovered characterization foliations for area-minimizing hypercones.
D. Ferus  and H. Karcher \cite{FK} showed, by constructing characterization foliation, that every cone over the minimal isoparametric hypersurface of an inhomogeneous isoparametric foliation on a sphere is area-minimizing.
G. Lawlor \cite{Law} 
completed the classification of all homogeneous area-minimizing hypercones.
Hence one can get a classification of all isoparametric homogeneous area-minimizing hypercones accordingly.
Actually, for each $C$ of these minimizing hypercones, L. Simon \cite{LS} gave a beautiful construction of minimal graph with tangent cone $C\times \mathbb R$ at infinity, 
thus creating a huge variety of solutions to the Bernstein problem.


\Section{Dirichlet problem for minimal surfaces of high codimensions}{Dirichlet problem for minimal surfaces of high codimensions}\label{S2}
Given an open bounded, strictly convex $\Omega\subset\mathbb R^{n+1}$ and $\phi:\p \Omega\rightarrow {\mathbb R}^{m+1}$,
           the \textbf{Dirichlet problem} (cf. \cite{j-s,b-d-m, de, m1,l-o}) 
           searches for weak solutions $f\in C^0(\bar \Omega)\cap Lip(\Omega)$ such that
             \begin{equation}\label{ms}
         \left\{\begin{array}{cc}
         \sum\limits_{i=1}^{n+1}\pf{x^i}(\sqrt{g}g^{ij})=0, & j=1,\cdots,n+1,\\
         \sum\limits_{i,j=1}^{n+1}\pf{x^i}(\sqrt{g}g^{ij}\pd{f^\a}{x^j})=0, & \alpha=1,\cdots,m+1,
         \end{array}
         \right.
         \end{equation}
 where $g_{ij}=\de_{ij}+\sum\limits_{\a=1}^{m+1}\pd{f^\a}{x^i}\pd{f^\a}{x^j}$, $(g^{ij})=(g_{ij})^{-1}$ and $g=\det(g_{ij})$, and further 
         $$
          f|_{\partial \Omega}=\phi.
          $$
       Note that $F(x)\mapsto (x, f(x))$ being harmonic (i.e., \eqref{ms}) is equivalent to its (or its graph) being minimal with respect to the induced metric from Euclidean space.
 %
When $m=0$,  \eqref{ms} can be reduced to the classical \eqref{DP1} to which many literatures were devoted as mentioned in \S \ref{S1}.
 
 In this section we shall talk about the case of $m\geq 1$.
 An astonishing pioneering work was done by Lawson and Osserman in \cite{l-o}, in which $\Omega$ is always assumed to be a unit disk $\mathbb D^{n+1}$. 
 In particular, they exhibited the following remarkable differences.
 \begin{itemize}
          \item [(1)] For $n=1$, $m\geq 1$, real analytic boundary data can be found
                     so that
                   there exist at least three different analytic solutions to the Dirichlet problem.
                   Moreover,
                   one of them has unstable minimal graph. 
          \item [(2)] For $n\geq 3$ and $n-1\geq m\geq 2$, the problem is in general not solvable.
                   A non-existence theorem is that, for each $C^2$ map $\eta:S^{n}\ra S^{m}$ that is not homotopic to zero under the dimension assumption,
                   there exists a positive constant $R_\eta$ depending only on $\eta$, such that
                   the problem is unsolvable for the boundary data $\phi=R\cdot \eta$, where
                   $R$ is a (vertical rescaling) constant no less than $R_\eta$.
          \item [(3)] For certain boundary data, there exists a Lipschitz solution to the Dirichlet problem which is not $C^1$.
\end{itemize}
        
        The ideas are briefly summarized as follows.
        
        (1) is based on a classical result by Rad\'o for $n=1$ case, 
        which
        says that every solution to the Plateau problem for boundary data given by a graph over boundary of a convex domain in some 2-dimensional plane
        has to be a graph over that domain.
        In fact, Lawson and Osserman were able to construct an action invariant boundary data of graph type for a $Z_4$-action in the total ambient Euclidean space $\mathbb R^{3+m}$ (for $m\geq 1$),
        such that under action of a generator of this $Z_4$-action each geometric solution to the Plateau problem (in the minimizing setting) cannot be fixed.
        Namely, one gains two distinct geometric solutions to that boundary, 
        and therefore, according to Rad\'o, two essentially different solutions to the corresponding Dirichlet problem.
        Then 
        by \cite{m-t} and \cite{Shi}
        there exists an unstable minimal solution of min-max type to the same boundary.
        Such boundary condition violates the uniqueness of solution and the minimizing property of solution graph.
        In particular, they constructed boundary supports at least three analytic solutions to the Dirichlet problem.
         It seems that more than 3 solutions may be created for certain boundary data, 
         if one considered symmetry by a discrete action of higher order group and actions of the entire group in some more subtle way.
        
        (2) is due to a nice special volume expression and the well-known density monotonicity for minimal varieties in Euclidean space.
       The proof of this meaningful result was achieved through a contradiction argument.
        Roughly speaking, the former can provide an upper bound for volume of graphs of solutions (as long as existed);
        while the latter guarantee a lower bound.
        Combined with the dimension assumption, these two bounds together lead to a contradiction
        when the rescaling factor becomes big.
        However, it is still completely mysterious and quite challenging to us how to figure out  the exact  maximal value of stretching factor with existence of solution(s).

        (3) is stimulated by (2).
       After establishing the non-existence result (2), Lawson-Osserman realized that, for a map satisfying both the dimension and homotopy conditions,
        if one rescaled the vertical stretching  factor by a tiny number,  then Dirichlet problem is solvable due to the Implicit Functional Theorem, e.g. see \cite{n};
        however if by a quite large number, then no Lipschitz solution can ever exist to the rescaled boundary functions.
        So a natural {\bf philosophy} by Lawson-Osserman states that
        there should exist $R_0$ such that the boundary condition $R_0\cdot\eta$ supports a singular solution.
        For first concrete examples of such kind,
        they considered the three noted Hopf maps between unit spheres.
        Expressed in complex coordinates
       $\eta(z_1,z_2)=(|z_1|^2-|z_2|^2,2z_1\bar z_2)$
       is the first.
       They looked for a minimal cone $C=C(\text{graph of }\phi)$ over graph $\phi=R_0\cdot\eta$. 
       Therefore, if existed, $C$ is also a graph with a link of ``spherical graph" type
       \begin{equation}
       \label{sg}
       L:=C\bigcap S^6=\left\{(\alpha x,\sqrt{1-\alpha^2}\eta(x):x\in S^3\right\}.
       \end{equation}
       Since a cone is minimal if and only if its link is a minimal variety in the unit sphere, it only needs to determine when $L$ is minimal.
       If one uses quaternions, then, isometrically up to a sign, $\eta(q)=qi\bar q$ for $q$ of unit length of $\mathbb H$ into pure imaginary part of $\mathbb H$,
       and $L$ can be viewed as an orbit through $((1,0,0,0), i)$ under action $Sp(1)\cong S^3$ with $q\cdot (\alpha x,\sqrt{1-\alpha^2}\eta(x))=(\alpha qx,\sqrt{1-\alpha^2}q\eta(x)\bar q)$.
       As a result, the orbit of maximal volume, corresponding to $\alpha={\frac{2}{3}}$, is minimal in $S^6$.
       Hence slope $R_0$ can take value $\frac{\sqrt{1-\alpha^2}}{\alpha}=\frac{\sqrt 5}{2}$.
       Similar procedures can be done for the other two Hopf maps.

        Inspired by the above, in recent joint work \cite{x-y-z0}, we attacked the question directly by generalizing \eqref{sg}.
        We introduced
        \begin{defi}\label{d1}
             A $C^2$ map  $\eta: S^{n}\rightarrow S^{m}$
           is called an Lawson-Osserman map {\bf (LOM)}
            if there exists $\theta\in (0,\frac{\pi}{2})$, 
            s.t.
        $F(x) :=(\cos\theta\cdot x,\sin\theta\cdot \eta(x))$
        gives a mininmal submanifold in $S^{m+n+1}$.
      The cone C(Image(F)) is  called associated Lawson-Osserman cone {\bf (LOC)}.
     \end{defi}
        
        \begin{rem}\label{r1}
                   For $\phi=\tan\theta\cdot \eta$, 
           $C(Graph(\phi))=C(Image(F))$ is a mininmal graph.
           So there is 
           a singular solution given by $f(x)=\begin{cases}
           |x|\cdot \tan\theta\cdot\eta(\frac{x}{|x|}),      & x\neq0;      \\
                   0,      & x=0.
\end{cases}$
        \end{rem}
        
       By Remark \ref{r1} it is clear that each Lawson-Osserman map induces a boundary function $\phi$ 
       which supports a cone-type singular solution.
       Then how many LOMs? In \cite{x-y-z0} we give a characterization.
       \begin{thm}\label{t1}
        A $C^2$ map $\eta:\ S^{n}
        {\rightarrow} S^{m}$ is LOM if and only if the followings hold
        for standard metrics $g_{m+n+1}, g_m, g_n$ of unit spheres
        \begin{equation}\label{eqc1}
          \begin{cases}
           \eta:(S^n,F^*g_{m+n+1})
           \rightarrow
           (S^m,g_m) 
           \text{ is harmonic};\\
           \sum_{i=1}^n\dfrac{1}{\cos^2\theta+\lambda_i^2\sin^2\theta}=n,
           \text{ where }\lambda_i^2  \text{ are diagonals of }\eta^*g_{m+n+1}  \text{ to } g_n.
\end{cases}
        \end{equation}
         \end{thm}

           In order to better understand the second condition in \eqref{eqc1},
           we put a strong restriction. 
        \begin{defi}  \label{d2}   
        $\eta$ is called an {\bf LOMSE}, if it is an LOM and in addition, for each $x\in S^n$, all nonzero singular values of $(\eta_*)_x$ are equal, i.e.,
$$\{\lambda_1,\cdots,\lambda_n\}=\{0,\lambda\}.$$
     \end{defi}
     \begin{rem}\label{r2}
              Let $p,\, n-p$ be the multiplicities for $\lambda$ and $0$.
              Then the second in \eqref{eqc1} becomes 
              $$\frac{n-p}{\cos^2\theta}+\frac{p}{\cos^2\theta+\lambda^2\sin^2\theta}=n.$$
              From the equality one can easily deduce that $p$ and $\lambda$ have to be independent of point $x$.
       \end{rem}
              So how many these LOMSEs? There turns out to be a constellation of uncountably many, even under the severe restriction!
              In \cite{x-y-z0} we derived a structure theorem.
            \begin{thm}\label{t2}
             $\eta$ is an LOMSE 
          if and only if
           $\eta=i\circ\pi$
           where $\pi$ is a Hopf fibration to $(\mathbb P^p,h)$ 
           and $i:(\mathbb P^p,\lambda^2h) 
           {\looparrowright} (S^m,g_m)$ is an isometric minimal immersion.
            \end{thm}
       \begin{rem}\label{r31}
       $\pi$ gives a countably many levels and 
       in most levels the moduli space of isometric minimal immersions from projective spaces into standard spheres 
       form (a sequence of) compact convex bodies of vector spaces of high dimensions.
      \cite{x-y-z0} perfectly embeds the relevant theory (see \cite{c-w,wa,oh,u,to,to2}) into the construction of LOMSEs.
         \end{rem}
       \begin{rem}\label{r32}
      In particular, using coordinates of ambient Euclidean spaces, 
       $\eta$ can be expressed as $(\eta_1,\cdots,\eta_{m+1})$.
       All $\eta_i$ are spherical harmonic polynomials sharing a common even degree $k$.
       Moreover $\lambda=\sqrt{\frac{k(k+n-1)}{p}}$.
       We call such an LOMSE of ${\bf (n,p,k)}$ type.
         \end{rem}
         
         Besides singular solutions, we cared about smooth solutions as well.
         By Morrey's  famous regularity result \cite{mo},
        a $C^1$ solution to \eqref{ms} is automatically $C^\omega$.
         In particular,
         a preferred variation of LOC associated to an LOM $\eta$ can be
$$
               M=M_{\rho,\eta}:=
               \{(rx,\rho(r)\eta(x)):x\in S^n, r\in 
               (0,\infty)
               \}
               \subset \mathbb R^{m+n+2}.
               $$
            Its being minimal is equivalent to two conditions (similar to that of \eqref{eqc1}, see \cite{x-y-z0} for details).
            When $\eta$ is an LOMSE,
             one of the conditions holds for free and the other gives the following. 
               
               \begin{thm}\label{t3}
For an LOMSE $\eta$, $M$ above is minimal if and only if
                                    \begin{equation}\label{ODE1}
\frac{\rho_{rr}}{1+\rho_r^2}+\frac{(n-p)\rho_r}{r}+\frac{p(\frac{\rho_r}{r}-\frac{\la^2\rho}{r^2})}{1+\frac{\la^2\rho^2}{r^2}}=0.
\end{equation}
               \end{thm}
               
               By introducing $\varphi:=\frac{\rho}{r}$ and  $t:=\log r$,
               \eqref{ODE1} transforms to
      \begin{equation}\label{ODE2}
      \left\{
      \begin{array}{ll}
\varphi_t=\psi,\\
\psi_t=-\psi-\Big[\big(n-p+\frac{p}{1+\la^2\varphi^2}\big)\psi+\big(n-p+\frac{(1-\la^2)p}{1+\la^2\varphi^2}\big)\varphi\Big]
\big[1+(\varphi+\psi)^2\big].
\end{array}
\right.
      \end{equation}
      {\ }\\
      This system is symmetric about the origin and owns exact 3 fixed points $(0,0),\, P(\varphi_0, 0)$ and $-P$,
       where $\varphi_0=\tan\theta$.
      Through linearization, it can be seen that the origin is always a saddle point
      and
      $P$ has two types:
           \begin{enumerate}
{
\item [(I)]
$P$ is a {stable center} when $(n,p,k)=(3,2,2), (5,4,2), (5,4,4)$ or $n\geq 7$;
}
{
\item [(II)] $P$ is a {stable spiral point} when $(n,p)=(3,2)$,
$k\geq 4$ or $(n,p)=(5,4)$, $k\geq 6$.
}
\end{enumerate} 

By very careful analysis including excluding limit circles,
there exists a special orbit emitting from the origin and approaching to $P$ for the system \eqref{ODE2}
and for $t\in (-\infty, +\infty)$. 
%
                    \begin{figure}[h]
                              \begin{minipage}[c]{0.4\textwidth}
                              \includegraphics[scale=0.46]{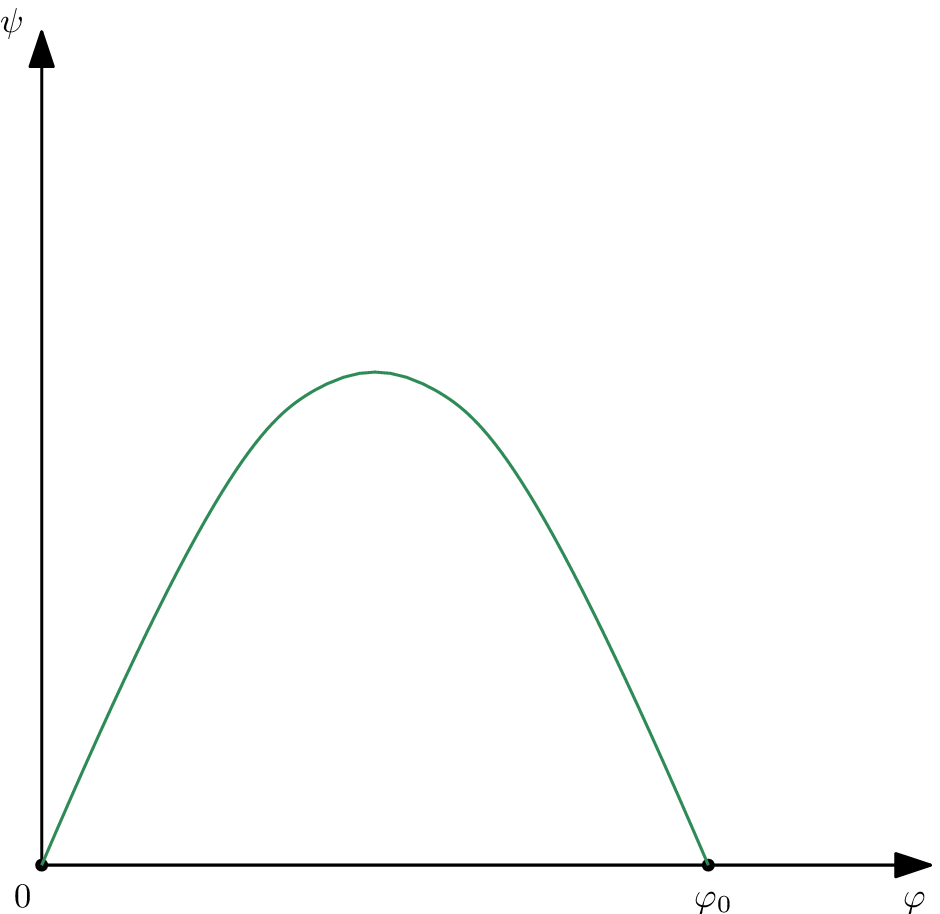}
                              \end{minipage}%
                          \begin{minipage}[c]{0.65\textwidth}
                           \includegraphics[scale=0.55]{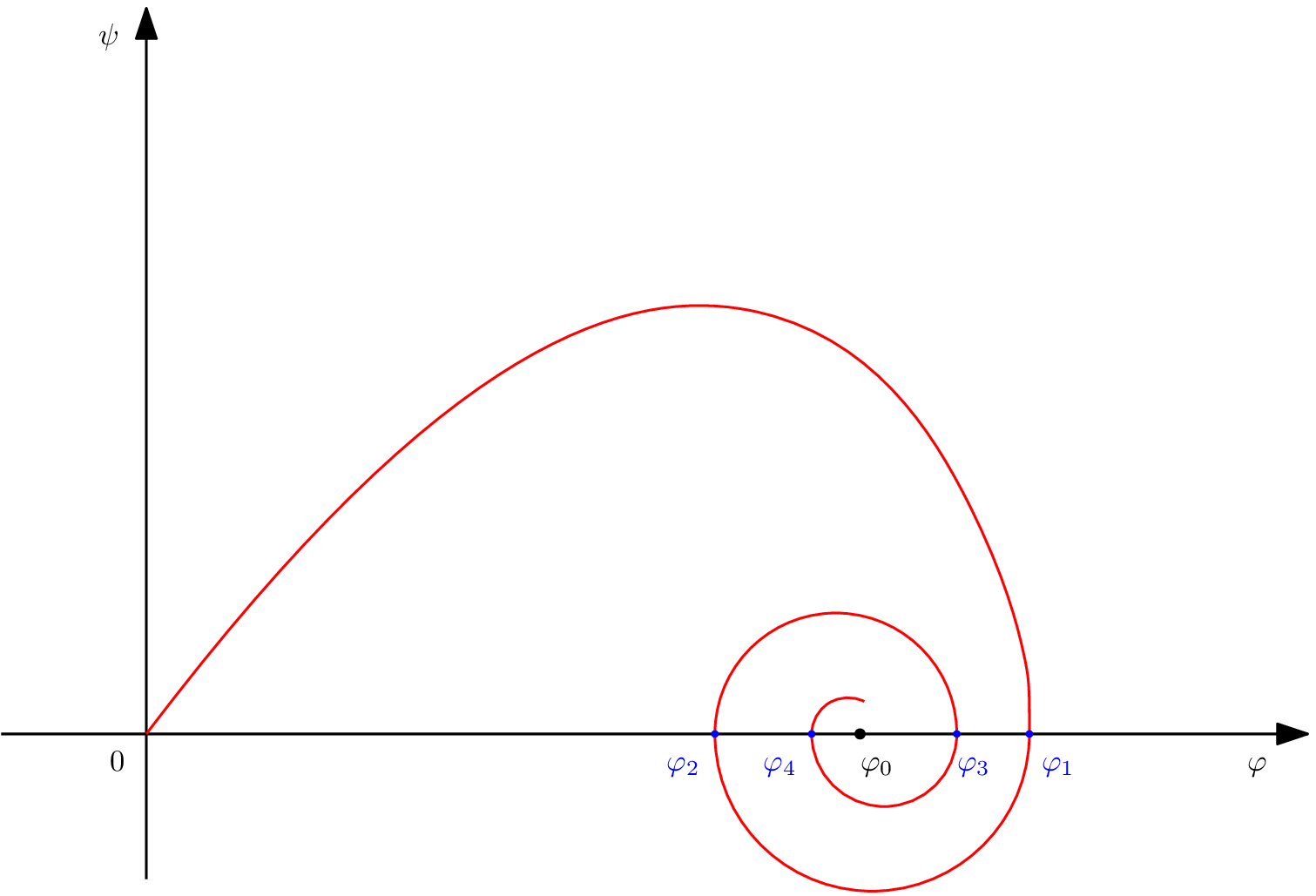}
                           \end{minipage}
                    \end{figure}

        Translated back to the $r\rho$-plane, the illustration graphs would be
        $$\begin{minipage}[c]{0.5\textwidth}
                              \includegraphics[scale=0.45]{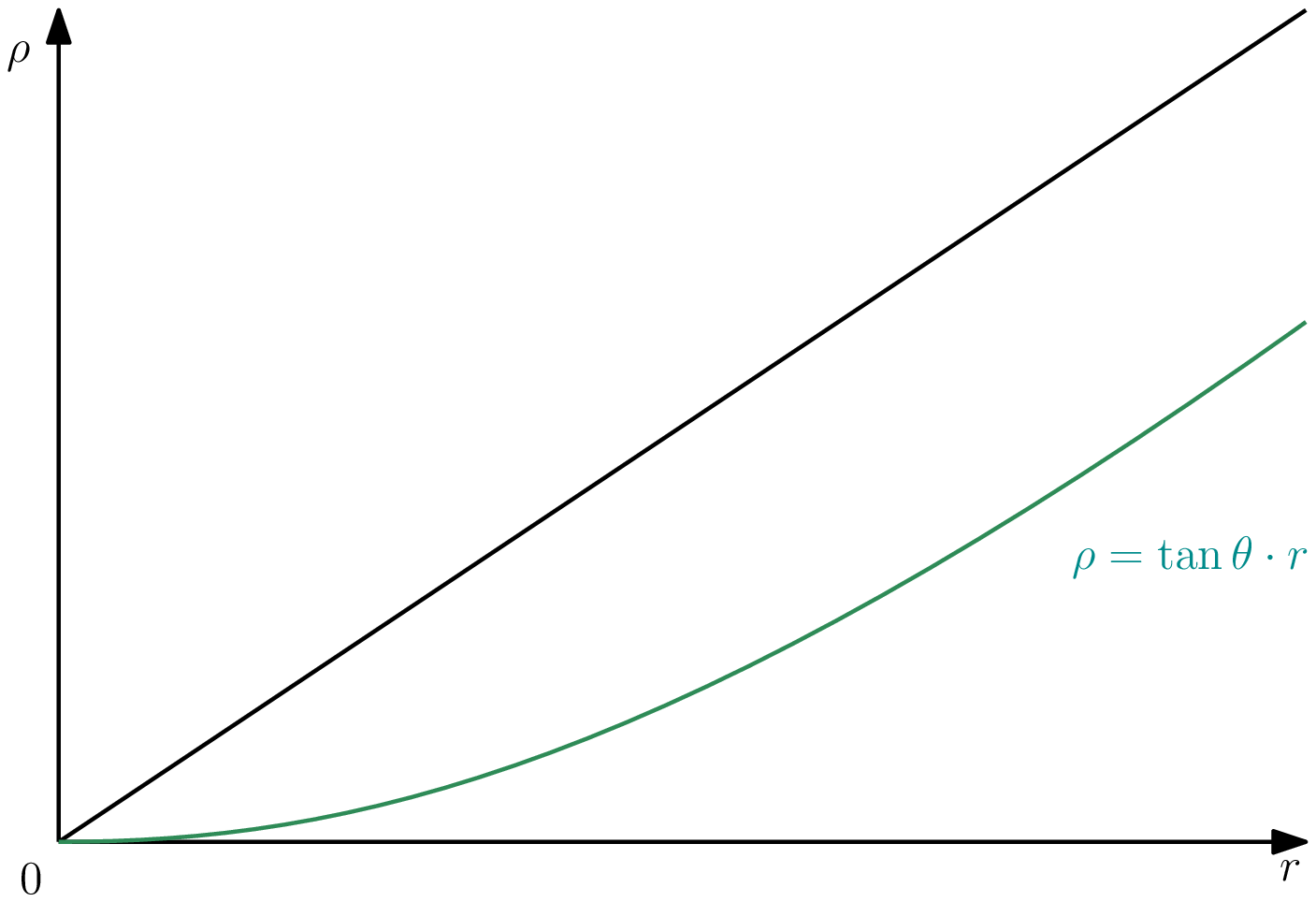}
                              \end{minipage}%
                          \begin{minipage}[c]{0.5\textwidth}
                           \includegraphics[scale=0.45]{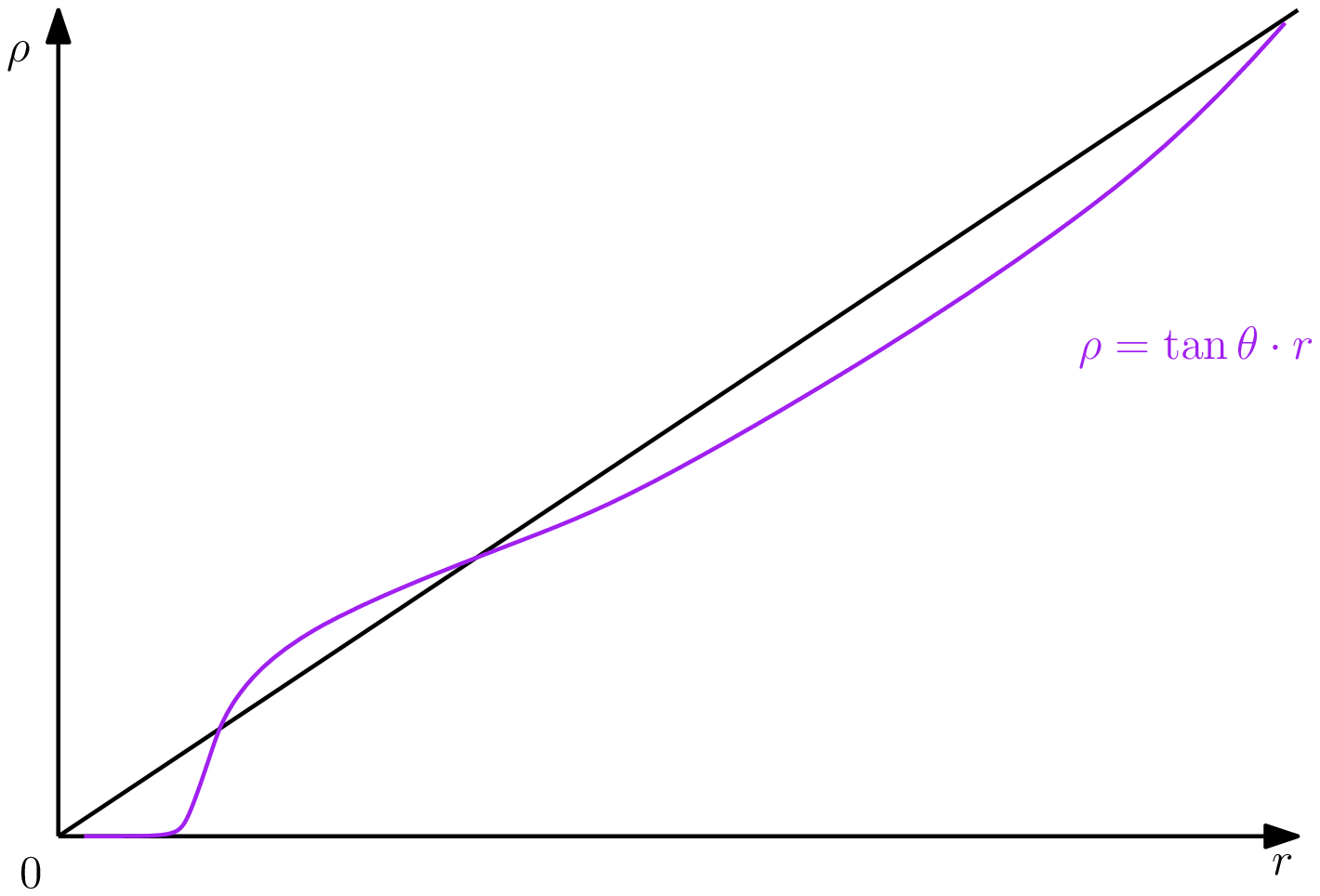}
                           \end{minipage}$$
                           
         Therefore, we got minimal graphs (other than cone) in $\mathbb R^{m+n+2}$ defined everywhere away from the origin of $\mathbb R^{n+1}$.
         Since $\frac{d\rho}{dr}=\varphi+\psi$,
         intense attentions were given to orbits emitting from the origin in $\varphi\psi$-plane.
          They produce minimal surfaces which are $C^1$ at $r=0$.
          So the natural $C^0$ extension is in fact $C^\omega$ according to Morrey's regularity result.
          This is the way how we constructed entire $C^\omega$ minimal graphs with LOCs as tangent cones at infinity.
          
          Type (II) contains interesting information. 
          Note that the fixed point orbit $P$ stands for the LOC and the ray in the $r\rho$-plane with constant slope $\varphi_0=\tan\theta$.
                   Since vertical line $\varphi=\varphi_0$ intersects the orbit infinitely many times,
         there are corresponding intersections of the solution curve and the LOC ray in $r\rho$-plane.
          Each intersection point gives us a minimal graph $G_i$ over a disk of radius $r_i$.
          Rescale $G_i$ by $\frac{1}{r_i}$ and denote new graphs by $\tilde G_i$.
          Then $\tilde G_i$ are mutually different minimal graphs over the unit disk $\mathbb D^{n+1}$
          with the same boundary $-$ graph of $\tan\theta\cdot\eta$ over unit sphere $S^n$.
          Hence we see that there exist boundary data which support infinitely many $C^\omega$ solutions and at least one singular solution to the Dirichlet problem!
          This extended Lawson-Osserman's non-uniqueness result (1) from finiteness to infiniteness.

          More can be read off. 
           Clearly, by density monotonicity, volumes of $\tilde G_i$ strictly increase to that of the truncated LOC.
          So none of the LOCs of Type (II) are area-minimizing.
          Actually, recently we
           showed 
           in a  joint work \cite{NZ}
           that LOC of Type (II) are even not stable.
            They bring unstable singular solutions to the Dirichlet problem (cf. (1) for Lawson-Osserman construction).
          Since the solution curve oscillates between rays of slopes $\varphi_1$ and $\varphi_2$,
          it is  not always the case that once we have singular solution for slope $\varphi_0$,
          then solutions suddenly vanish  immediately  for $\varphi>\varphi_0$.
                   It is a question for what kind of $\varphi$ outside $[0,\varphi_1]$ the problem can be solved?
          Union of the set of such value and $[0,\varphi_1]$  is called {\bf slope-existence range} of $\eta$ for the Dirichlet problem.
          To extend $[0,\varphi_1]$,
          maybe the first difficulty is to figure out whether the orbit between the origin and point $(\varphi_1,0)$ gives a stable compact minimal graph.

         In the opposite direction on non-existence,
         a recent preprint \cite{z0} confirmed that the slope-existence range 
         should usually be contained in a compact set of $\mathbb R_{\geq 0}$.
        More precisely, we prove
        \begin{thm}
        For 
       every
        LOMSE $\eta$ of         
         {either Type (I) or Type (II)},
         there exists positive constant $R_\eta$ such that
         when constant $R\geq R_\eta$,
         the Dirichlet problem has no solutions for $\phi=R\cdot \eta$.
        \end{thm}
        

\section{Minimal cones}\label{S3}
As briefly mentioned before, it is useful to see if a local structure is stable or not for observation, and
it is also quite important to know structures of minimizing currents.
Minimal cones are infinitesimal structures of minimal varieties, while minimizing cones are infinitesimal structures of minimizing currents.
Both determine, in some sense, local diversities of certain geometric objects.

In fact we naturally encountered many examples of minimal cones.
For example, Lawson-Osserman \cite{l-o} constructed three minimal cones for singular solutions to the Dirichlet problem.
The first cone was shown to be coassociative in $\mathbb R^7$ and hence area-minimizing by the fundamental theorem of calibrated geometries
 in the milestone paper \cite{h-l}.
 However, it was unknown for 40 years if the other two are minimizing or not.
 In our recent joint work \cite{x-y-z}  we proved that all LOCs of $(n, p, 2)$ type (for which case moduli space of $i$ to each Laplacian eigenvalue is a single point)
 are area-minimizing.
 Since the other two original Lawson-Osserman cones are of $(7,4,2)$ type and $(15,8,2)$ type respectively, 
 so the long-standing question got settled.

 Area-minimizing cones of $(n,p,2)$ type are all homeomorphic to Euclidean spaces. 
 For other kind of area-minimizing cones, we considered those associated to isoparametric foliations of unit spheres.
 There are two natural classes of minimal surfaces $-$ minimal isoparametric hypersurfaces and focal submanifolds.
By virtue of a successful combination of Lawlor's curvature criterion and beautiful structure of isoparametric foliations,
we were able  in \cite{TZ}  to show that, except in low dimensions, cones overs the ``minimal products" (defined therein) 
among these two classes are area-minimizing.
These provide a large number of new area-minimizing cones with various links of rich complexities.
Note that none of them can be split as product of (area-minimizing) cones of lower dimensions.
It is currently unknown to the author whether minimal products of links of general area-minimizing cones  can always span an area-minimizing cone.

In \cite{z} we considered a realization problem, first attacked by N. Smale \cite{NS, NS2} in later 1990s.
\begin{quote}
{\it Can any area-minimizing cone be realized as a tangent cone at a point
of some homologically area-minimizing  {\tt compact} singular submanifold?}
\end{quote}
N. Smale constructed first such examples by applying many tools in geometric analysis and geometric measure theories in \cite{NS},
while ours seems a bit simpler through the theory of calibrations with necessary understandings on Lawlor's work \cite{Law}. 
We showed
\begin{thm}\label{t4}
Every oriented area-minimizing cone in \cite{Law} can be realized to the above question.
\end{thm}
\begin{rem}\label{r4}
Prototypes can be all the newly-discovered oriented area-minimizing cones in \cites{TZ, x-y-z} and
all Cheng's examples of homogeneous area-minimizing cones of codimension $2$ in \cite{Ch}
(e.g. minimal cones over
                        $\text{U}(7)/\text{U}(1)\times \text{SU}(2)^3$ in $\mathbb R^{42}$,
                        $\text{Sp}(n)\times \text{Sp}(3)/\text{Sp}(1)^3\times \text{Sp}(n-3)$ in $\mathbb R^{12n}$ for $n\geq 4$,
                        and $\text{Sp}(4)/\text{Sp}(1)^4$ in $\mathbb R^{27}$)
                        via a variation of  our arguments in \cite{z}.
\end{rem}

All the above cones have smooth links. It could be highly useful if one can derive an effective way to study cones with non-smooth links.

As for stability and instability, in \cite{NZ} we borrowed ideas 
 \cite{br, h-l, l} for orbit space.
          We focused on a preferred subspace associated to given LOMSE
          and its quotient space.
          With a canonical metric
          $
          \sigma_0^2\cdot 
               \left[
               \left(
               r^2+\lambda^2\rho^2
               \right)^p
               \cdot
               r^{2(n-p)}
               \right]
               \cdot [dr^2+d\rho^2]
               $
               where $\sigma_0$ is the volume of $n$-dimensional unit sphere,
          the length of any curve in the quotient space 
          equals the volume of corresponding submanifold in $\mathbb R^{m+n+2}$.
          Hence, the infinitely many $C^\omega$ solution curves in the $r\rho$-plane for Type (II) in \S \ref{S2} determine geodesics 
          connecting $Q:=(1,\tan\theta)$ and the origin in the quotient space.
                    \begin{figure}[h]
          \includegraphics[scale=0.55]{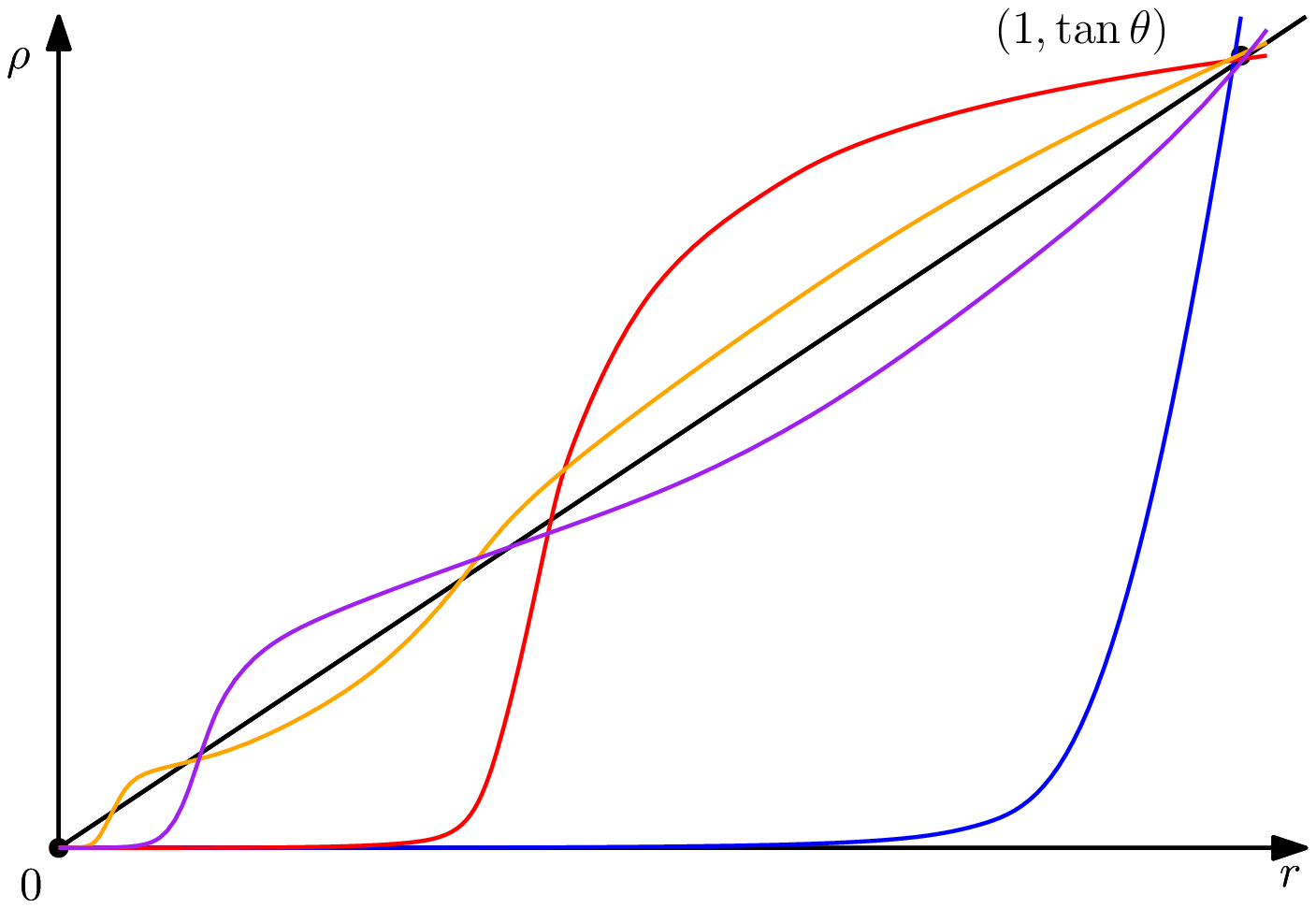}
          \label{AS}
         \end{figure}
         \\
         {\ }
         \\
We showed
          \begin{thm}\label{t5}
  The line segment $\overline{0Q}$ is stable for Type \text{(I)} and unstable for Type {(II)}.
  \end{thm}
   \begin{rem}\label{r5}
  In fact, $\overline{0Q}$ is minimizing for Type \text{(I)}.
 However, the difficulty is whether it is possible to lift the stability or even the area-minimality property back to LOCs in $\mathbb R^{m+n+2}$.
  \end{rem}

\section{Open questions}\label{S4}
Besides several open questions in previous sections, we want to emphasize a few more in this section.

1. Systematic study about LOMs beyond LOMSEs.
It still remains unclear how to construct LOMs for other type distribution of critical values, for instance 
$\{0,\, \lambda_1,\, \lambda_2\}$ or $\{\lambda_1,\, \lambda_2\}$ where $\lambda_1$ and $\lambda_2$ are different positive numbers.
These may involve more complicated dynamic systems and perhaps chaos phenomena.

2. It seems unknown in general if the cones over image of $i$ itself in Remark \ref{r31} is area-minimizing or not.
This would need certain systematic understandings about second fundamental form of $i$.
It can also help us a lot for a complete classification of which Lawson-Osserman cones associated to LOMSEs
are area-minimizing. 

3. How about situation for  the finite left cases in \cite{TZ}.
The most famous one may be the cone over the image of Veronese map in $S^4$, a focal submanifold for the isoparametric foliation with $g=3$ and $m=1$.
 It is still open whether the cone over the image, a minimal embedded $\mathbb RP^2$ of constant curvature, is a minimizing current mod 2 (see \cite{Zm}).

{\ }
 
 \section*{Aknowlegement}
 The author would like to thank MPIM at Bonn for warm hospitality.
 This work was sponsored in part by {the}
NSFC (Grant No. {11601071})
 and a Start-up Research Fund from Tongji University. 

{\ }

\end{document}